\documentclass[11pt,reqno]{amsart}
\usepackage[french,english]{babel}
\usepackage[utf8]{inputenc}
\usepackage[T1]{fontenc}            
\usepackage{amsmath}
\usepackage{graphicx}
\usepackage{amssymb}
\usepackage{amsthm}
\usepackage{comment}
\usepackage{pgf, tikz}
\usepackage{amsfonts,mathrsfs}
\usepackage{thmtools}
\usepackage{geometry}
\usepackage{stmaryrd}
\usepackage{color}
\usepackage{hyperref}
\usepackage{bm}
\hypersetup{colorlinks = true , linkcolor=blue, citecolor=blue}
\numberwithin{equation}{section}
\usepackage{comment}
\theoremstyle{definition}

\declaretheorem[numberwithin=,title=Proposition]{Prop}

\declaretheorem[numberwithin=,title=Theorem]{Thm}

\declaretheorem[numberwithin=,title=Remark]{Rq}

\newtheorem{Assumption}{Assumption}

\newcommand{\FF}{\mathcal{F}}
\newcommand{\E}{\mathbb{E}}

\newcommand{\R}{\ensuremath{\mathbb{R}}\xspace}

\newcommand{\CC}{\mathcal{C}}

\newcommand{\1}{\mathbf{1}}
\newcommand{\MM}{\mathcal{M}}
\newcommand{\NN}{\mathcal{N}}
\newcommand{\NNN}{\widetilde{\mathcal{N}}}
\renewcommand{\P}{\mathbb{P}}

\newcommand{\Pc}{\mathcal{P}}
\newcommand{\PP}{\mathcal{P}_{\beta}(\mathbb{R}^d)}
\newcommand{\PPP}{\mathcal{P}(\mathbb{R}^d)}

\newcommand{\muu}{\overline{\mu}}

\begin{document}
	
	\title[Lévy-driven McKean-Vlasov SDEs]{Well-posedness and propagation of chaos for Lévy-driven McKean-Vlasov SDEs under Lipschitz assumptions}

	\author{Thomas Cavallazzi}
	\address{Université Paris-Saclay, CentraleSupélec, MICS and CNRS FR-3487, France.}
	\email{thomas.cavallazzi@centralesupelec.fr}
	\date{January 19, 2025}

	\begin{abstract}The first goal of this note is to prove the strong well-posedness of McKean-Vlasov SDEs driven by Lévy processes on $\R^d$ having a finite moment of order $\beta \in [1,2]$ and under standard Lipschitz assumptions on the coefficients. Then, we prove a quantitative propagation of chaos result at the level of paths for the associated interacting particle system, with constant diffusion coefficient. Finally, we improve the rates of convergence obtained for linear interactions with respect to the measure and when the noise is a $\alpha$-stable process with $\alpha \in (1,2)$, for which we have $\beta < \alpha$.\end{abstract}
	
	\maketitle

    \section{Introduction and results}

Let us fix $(\Omega,\FF,(\FF_t)_{t\geq 0},\P)$ a filtered probability space and $\NN$ a Poisson random measure on $\R^+ \times \R^d\backslash\{0\}$ with intensity $dt \otimes \nu$, where $\nu$ is a Lévy measure, i.e. $$\nu(\{0\})=0 \quad \text{and}\quad \int_{\R^d} 1\wedge |z|^2\, d\nu(z) < +\infty,$$ where $a\wedge b$ denotes the minimum between to real numbers $a$ and $b$. We denote by $\NNN(dt,dz):= \NN(dt,dz) - dt\otimes d\nu(z)$ the associated compensated Poisson measure. We consider  $Z=(Z_t)_{t\geq 0}$ a Lévy process on $\R^d$ written, for all $t \geq 0$, as \begin{equation*}
	Z_t = \int_0^t \int_{B_1} z \, \NNN(ds,dz) + \int_0^t \int_{B_1^c} z \, \NN(ds,dz),\end{equation*}  where $B_1$ is the open ball of $\R^d$ centered at $0$ and of radius $1$ and $B_1^c$ is its complementary in $\R^d$.\\

We assume that there exists $\beta \in [1,2]$ such that the Lévy measure $\nu$ satisfies 
\begin{equation*}
	\int_{B_1^c} |z|^\beta \, d\nu(z) < + \infty.\end{equation*}
This is equivalent to assuming that for any $t \in \R^+$, $Z_t$ has a finite moment of order $\beta$ by \cite[Theorem $25.3$]{SatoLevyprocessesinfinitely1999}. Let us denote by $\PP$ the space of probability measures on $\R^d$ having a finite moment of order $\beta$, which is endowed with the Wasserstein metric $W_{\beta}$. We are interested in the well-posedness of the following Lévy-driven McKean-Vlasov SDE \begin{equation}\label{Aedsmckv}
	\left\{  \begin{array}{lll}
		&dX_t = b_t(X_t,\mu_t) \,dt + \sigma_t(X_{t^-},\mu_t)\, dZ_t,\quad t\in[0,T], \\ &\mu_t := [X_t],\\ &X_0 = \xi \in L^\beta(\Omega,\FF_0;\R^d),
	\end{array}\right.
\end{equation}
where $T$ is a fixed finite horizon of time, $[X_t]$ denotes the distribution of $X_t$ and $b: [0,T]\times \R^d \times \PP \to \R^d$ and $\sigma : [0,T]\times \R^d \times \PP \to \mathcal{M}_d(\R)$ are measurable maps, $\mathcal{M}_d(\R)$ being the space of matrices of size $d\times d$ on $\R$. The first motivation to study \eqref{Aedsmckv} lies into its connection with the following mean-field interacting particle system \begin{equation}\label{particle_sytm}
	\left\{  \begin{array}{lll}
		&dX^{i,N}_t = b_t(X^{i,N}_t,\muu^N_t)\,dt + \sigma_t(X^{i,N}_{t^-},\muu^N_t)\,dZ^i_t, \quad t \in [0,T],\quad i \in \{1,\dots,N\}, \\ &\muu^N_t := \frac{1}{N} \sum\limits_{j=1}^N \delta_{X^{j,N}_t},\\ &X^{i,N}_0 = \xi^i,
	\end{array}\right.\end{equation}
where $(Z^i,\xi^i)_{ i \geq 1}$ are i.i.d.\ with same distribution as $(Z,\xi)$. The link between \eqref{Aedsmckv} and \eqref{particle_sytm} is that for any $k \geq 1$, the dynamics of $k$ particles is expected to be described by $k$ independent copies of \eqref{Aedsmckv} when the total number of particles $N$ tends to infinity. This is the so-called propagation of chaos phenomenon. It was originally studied by McKean \cite{McKeanPropchaos67} and then investigated by Sznitman \cite{SznitmanTopicspropagationchaos1991} when $Z$ is a Brownian motion. For a detailed review on the topic of propagation of chaos, we refer the reader to \cite{chaintron:hal-03585065,chaintron:hal-03585067}. Let us mention some applications of these mean-field systems. They naturally appear for example in physics (kinetic theory), in biology to describe the motion of a cell population, in neuroscience to model the interactions between neurons and also in the Mean-Field Games theory. In particular these systems driven by Lévy processes are used to model some physical systems (Lévy flights and anomalous diffusions, see \cite{Mann} and \cite{Jourdainfractional} for example). \\

We are going to work under the following Lipschitz assumptions.

\begin{Assumption}\label{Assumpt_lip}
	There exists a constant $C>0$ such that for all $t \in [0,T],$ $x,y \in \R^d$ and $\mu,\nu \in \Pc_{\beta}(\R^d),$ we have \begin{equation}\label{assumpt_lip_eq}
		|b_t(x,\mu) - b_t(y,\nu)| + |\sigma_t(x,\mu) - \sigma_t(y,\nu)| \leq C(|x-y| + W_{\beta}(\mu,\nu)),
	\end{equation}
	and 
	\begin{equation*}
		|b_t(x,\mu)| + |\sigma_t(x,\mu)| \leq C(1+|x| + M_{\beta}(\mu)),
	\end{equation*}
	where $M_{\beta}(\mu) = \left(\int_{\R^d} |x|^\beta \, d\mu(x)\right)^{\frac{1}{\beta}}$ for $\mu \in \PP$.
\end{Assumption}

\subsection{Well-posedness of the McKean-Vlasov SDE \eqref{Aedsmckv}}

\begin{Thm}\label{ThmexistenceMcKV1}
	Let us recall that $\beta \in [1,2]$. Under Assumption \ref{Assumpt_lip}, there exists a unique strong solution $(X_t)_{t \in [0,T]}$ to \eqref{Aedsmckv} for all initial datum $\xi \in L^{\beta}(\Omega,\FF_0;\R^d).$ Moreover, the flow of marginal distributions $(\mu_t)_{t \in [0,T]}$  belongs to $\CC^0([0,T];\Pc_{\beta}(\R^d))$ and we have \begin{equation}\label{moment_estimate} \E \sup_{t\leq T} |X_t|^{\beta}  <+ \infty.\end{equation}
\end{Thm}

The proof of this theorem is based on the Banach fixed point theorem on the complete space $\CC^0([0,T],\mathcal{P}_{\beta}(\R^d)).$ To prove the contraction property, we need to prove that the map which associates to a flow $(\mu_t)_{t \in [0,T]} \in \CC^0([0,T],\mathcal{P}_{\beta}(\R^d))$ the flow of marginal distributions of the solution to \eqref{Aedsmckv}, where the measure argument of the coefficients is equal to $(\mu_t)_{t \in [0,T]}$, is a contraction. This is done by separating the small and large jumps of the noise by conditioning with respect to the large jumps, as done in \cite{fournier_moment}. Indeed, the small jumps cannot be treated in $L^{\beta}$ in general due to the integrability of the Lévy measure near $0$. This is the main technical difficulty in comparison with the Brownian case. 

\begin{Rq} We can easily add a term of the form $ (Bt + \Sigma W_t)_{t\geq 0}$ to $Z$, where $B \in \R^d$, $\Sigma\in \MM_d(\R)$ is a symmetric positive semidefinite matrix of size $d\times d$ and $W$ is a standard Brownian motion on $\R^d$. Using the Lévy-Itô decomposition given in \cite[Theorem $2.4.16$]{ApplebaumLevyprocessesstochastic2009}, we can thus consider a general Lévy process $Z$ having a finite moment of order $\beta\in[1,2]$.

\end{Rq}

Let us compare our result with the existing literature. When $\beta =2$, the well-posedness of \eqref{Aedsmckv} was proved by Jourdain, Méléard and Woyczynski \cite{JourdainNonlinearSDEsdriven2007}. In this work, the weak existence is also proved when $\beta =0$ through the relative nonlinear martingale problem. However, uniqueness is not shown when $\beta =0$. When $\beta=1$, a result similar to Theorem \ref{ThmexistenceMcKV1} is proved by Graham in \cite[Theorem $2.2$]{GrahamNonlineardiffusionjumps1992}. The main differences are the following. Firstly, in \cite{GrahamNonlineardiffusionjumps1992}, there is no integral with respect to the compensated Poisson random measure $\NNN$ in the definition of $Z$. Secondly, in the case where the drift $b$ is unbounded, it is supposed in \cite{GrahamNonlineardiffusionjumps1992} that $X_0$ has a finite moment of order $2$, which is not the case in Theorem \ref{ThmexistenceMcKV1}. Thirdly, it is also assumed that, keeping our notations, there exists $C>0$ such that for all $t \in [0,T]$, $x \in \R^d$ and $\mu \in \mathcal{P}_1(\R^d)$, we have \begin{equation}\label{_hyp_graham}\left\vert\int_{B_1^c}\sigma_t(x,\mu) z \, d\nu(z) \right\vert^2  + \int_{B_1^c} |\sigma_t(x,\mu)z|^2 \, d\nu(z) \leq C (1+|x|^2).\end{equation}
This assumption is more restrictive than those made in Theorem \ref{ThmexistenceMcKV1}. For example, when $\sigma = \text{id}$, it implies that $$\int_{B_1^c}|z|^2 \, d\nu(z) < + \infty.$$ But this is equivalent to the fact that for any $t \in \R^+$, $Z_t$ has a finite moment of order $2$, which is not supposed in Theorem \ref{ThmexistenceMcKV1} since $\beta \in [1,2]$. Moreover, the assumptions of Theorem \ref{ThmexistenceMcKV1} allow to consider unbounded coefficients with respect to the measure, which is not compatible with \eqref{_hyp_graham}, for example in dimension $d=1$. Our result is more general.\\

 In the non-degenerate case, i.e.\ when $\sigma$ is uniformly elliptic, we refer to \cite{frikha_menozzi_konakov}. In this work, Frikha, Menozzi and Konakov prove the well-posedness of \eqref{Aedsmckv} under Hölder assumptions on the coefficients with respect to both space and measure variables. Of course, this result can be applied to Lipschitz continuous coefficients but in Theorem \ref{ThmexistenceMcKV1}, we do not assume that the diffusion coefficient $\sigma$ is uniformly elliptic. Moreover, another assumption made in \cite{frikha_menozzi_konakov} is that for all $(t,x) \in [0,T]\times \R^d$, the maps $\mu\in \PPP \mapsto b_t(x,\mu)$ and $\mu \in \PPP \mapsto \sigma_t(x,\mu)$ have bounded linear derivatives, where $\PPP$ is the space of probability measures on $\R^d$. Note that, at least when the coefficients depend linearly on the measure, this assumption implies the boundedness of the coefficients with respect to the measure variable. This is not the case in the present work.

\begin{Rq}\label{_rq_beta_petit}
	Notice that when $\beta \in (0,1),$ the uniqueness result of Theorem \ref{ThmexistenceMcKV1} is false without a non-degeneracy assumption on the diffusion coefficient $\sigma.$ Let us give a simple counterexample by setting, for $t \in [0,T],$ $x \in \R^d$ and $\mu \in \PP$ $$ b_t(x,\mu) := \int_{\R^d} |x|^\beta \, d\mu(x), \quad \sigma_t(x,\mu) := 0,\quad \text{and} \quad\xi := 0.$$ Assumption \ref{Assumpt_lip} is clearly satisfied. Moreover, the solution to the corresponding McKean-Vlasov SDE is deterministic since there is no noise and the initial distribution is deterministic. We easily remark that the problem is equivalent to solve the ordinary differential equation $$\begin{cases} &y'(t) = |y(t)|^{\beta}, \quad t \in [0,T],\\ &y(0)=0.\end{cases}$$ It is well-known that there exist several solutions to this problem. However, under Assumption \ref{Assumpt_lip}, there exists at least one strong solution to the McKean-Vlasov SDE \eqref{Aedsmckv}. We refer to Appendix \ref{appendix_existence} for a proof of this result.
\end{Rq}

\subsection{Propagation of chaos for the interacting particle system \eqref{particle_sytm}}

We now focus on the propagation of chaos for the interacting particle system \eqref{particle_sytm}. Under Assumption \ref{Assumpt_lip}, the SDE \eqref{particle_sytm} admits a unique strong solution by \cite[Theorem $6.2.9$]{ApplebaumLevyprocessesstochastic2009}. Propagation of chaos can be understood in the weak sense, i.e.\ in distribution through the convergence of the empirical measure $\muu^N$, or in the strong sense, i.e.\ at the level of paths. Our aim is to prove quantitative strong propagation of chaos. Let us introduce the i.i.d.\ copies of  the limiting McKean-Vlasov SDE \eqref{Aedsmckv}, which are denoted by $(X^{i,\infty})_{i \geq 1}$, where the initial data and the noises are respectively $(\xi^i)_{i\geq 1}$ and $(Z^i)_{i\geq 1}$.

\begin{Thm} \label{thm_poc_lip}  Assume that $\beta \in (1,2]$ and that $\sigma = \text{Id}$. 
	 \begin{enumerate} \item Assume moreover that Assumption \ref{Assumpt_lip} holds true with $W_1$ instead of $W_\beta$ in the Lipschitz control \eqref{assumpt_lip_eq}. Then, there exists a constant $C>0$, independent of $d$ and $N$, such that for all $N \geq 1$
	\begin{equation}\label{eq1_rq_exemple}
		\sup_{i\leq N}\E  \sup_{t \leq T} |X^{i,N}_t - X^{i,\infty}_t| \leq C \begin{cases}
			N^{\frac{1}{\beta}-1}, \, &\text{if}\quad d=1,2\quad  \text{or}\quad d\geq3\,\, \text{and}\,\,\beta < \frac{d}{d-1},\\ 
			N^{-\frac1d}, \, &\text{if}\quad d\geq 3\,\, \text{and}\,\, \beta > \frac{d}{d-1},
		\end{cases}
	\end{equation}
	and \begin{equation}\label{eq2_rq_exemple}
		\sup_{t \in [0,T]}\E W_1(\muu^N_t,\mu_t) \leq C \begin{cases}
			N^{\frac{1}{\beta}-1}, \, &\text{if}\quad d=1,2\quad  \text{or}\quad d\geq3\,\, \text{and}\,\,\beta < \frac{d}{d-1},\\ 
			N^{-\frac1d}, \, &\text{if}\quad d\geq 3\,\, \text{and}\,\, \beta > \frac{d}{d-1}.
		\end{cases}
	\end{equation}

\item In the case where $b$ is of the form $b_t(x,\mu) = \int_{\R^d} \tilde{b}(x,y) \, d\mu(y),$ where $\tilde{b}$ is Lipschitz continuous, there exists a positive constant $C$ such that 
\begin{equation}\label{eq3_rq_exemple}
	\sup_{i\leq N}\E  \sup_{t \leq T} |X^{i,N}_t - X^{i,\infty}_t| \leq C N^{\frac1\beta -1}.
\end{equation}
Let us assume moreover that $Z=(Z_t)_{t \geq 0}$ is a $\alpha$-stable process with $\alpha \in (1,2)$ and that $\xi \in L^{\alpha}(\Omega,\FF_0,\P)$. Then, there exists a positive constant $C$ such that for all $N \geq 1$ \begin{equation}\label{eq1_thm_poc_exemple}
	\sup_{i\leq N}\E \sup_{t \in [0,T]}|X^{i,N}_t - X^{i,\infty}_t| \leq C N^{\frac1\alpha - 1} (\ln(N))^{\frac{1}{\alpha}}.
\end{equation}
\end{enumerate}
\end{Thm}

\begin{Rq}Note that by exchangeability, the supremum in \eqref{eq1_rq_exemple}, \eqref{eq3_rq_exemple} and \eqref{eq1_thm_poc_exemple} is free since the expectation does not depends on $i$. Let us also point out that the rate in \eqref{eq3_rq_exemple} is better than the one given in \eqref{eq1_rq_exemple}, where no structural assumption is made on the drift. Moreover, in the $\alpha$-stable case, if we directly apply \eqref{eq3_rq_exemple}, we have to choose $\beta \in [1, \alpha)$ since the moment of order $\alpha$ of $Z_t$ is infinite (see \cite[Theorem $25.3$]{SatoLevyprocessesinfinitely1999}). This leads to a bigger rate of propagation of chaos in \eqref{eq1_thm_poc_exemple} than in \eqref{eq3_rq_exemple}.
\end{Rq}

The proof of \eqref{eq1_rq_exemple} and \eqref{eq2_rq_exemple} is standard and similar to the Brownian case but we write it for the sake of completeness. It relies on the estimation of $\sup_{i \leq N}	\E  \sup_{t \leq T} \vert X^{i,N}_{t} - X^{i,\infty}_{t}|$ using Grönwall's inequality and the rate of convergence of the empirical measure associated with the i.i.d.\ random variables $(X^{i,\infty}_t)_N$ proved in \cite{fournier:hal-00915365}. The proof of \eqref{eq3_rq_exemple} is inspired by \cite{SznitmanTopicspropagationchaos1991} and relies also on the $\beta$-th order Marcinkiewicz-Zygmund inequality (see $(1.4)$ in \cite{ChenSunginequalitysum}).
The proof of \eqref{eq1_thm_poc_exemple} relies on the technical trick of truncating the jumps of size larger than the number of particles $N$ from all the noises of $X^{i,N}$ and $X^{i,\infty}$. We then control the corresponding error and use the same arguments as in the proof of \eqref{eq3_rq_exemple}.

\begin{Rq}
The method used in the proof of Theorem \ref{thm_poc_lip} cannot be applied to prove quantitative strong propagation of chaos with a general non-constant diffusion coefficient $\sigma$ under Assumption \ref{Assumpt_lip}. This is a bit disapointing since in the Brownian case, adding a Lipschitz diffusion coefficient is not a difficulty. The difference in the Lévy-driven case is that we have to treat differently the small and large jumps due to the integrability property of the Lévy measure near $0$. One could try to use the same strategy as in the proof of Theorem \ref{ThmexistenceMcKV1} but the separation of the small and large jumps and the conditioning arguments in this proof are much more difficult to handle when applied to the particle system. Thus, it remains, to the best of our knowledge, an open problem.
\end{Rq}

We now compare our result with the existing literature. In \cite{GRAHAM199269}, Graham proves non-quantitative weak propagation of chaos, i.e.\ without rate of convergence, under Lipschitz assumptions, similar as Assumption \ref{Assumpt_lip}, for an interacting particle system driven by a Poisson random measure and its compensated
measure. It is supposed that the Poisson random measure is associated with a Poisson process
having a finite moment of order $1$. In this work, $\sigma$ is not assumed to be the identity matrix, but the result is only qualitative. In the framework of SDEs driven only by Poisson random measures, let us mention \cite{DaiPra_PoC} and \cite{Erny_PoC} for quantitative strong propagation of chaos results with $\sigma \neq \text{Id}$. In our framework, the main difficulty comes from the compensated Poisson random measure associated to the small jumps of the Lévy process. This term can be easily treated in the $L^2$ framework using the $L^2$ isometry but it is much more difficult to treat it in the $L^\beta$ framework, with $\beta <2$, which is the case here. For a general Lévy noise having a finite moment of order $2$ and still in the Lipschitz framework, we refer to Jourdain et al. \cite{JourdainNonlinearSDEsdriven2007} and Neelima et al. \cite{NeelimaWellposednesstamedEuler2020}, where quantitative strong propagation of chaos is proved in $L^2$. The assumptions are a bit weaker in \cite{NeelimaWellposednesstamedEuler2020}, as well as compared to those of Theorem \ref{thm_poc_lip}, since they allow super-linear growth in the state
variable. Their rates of convergence also come from \cite{fournier:hal-00915365} and are thus better than those of Theorem \ref{thm_poc_lip} since they can work in $L^2$. In these works, $\sigma$ is not assumed to be the identity matrix but the $L^2$ framework allows to do the same kind of computations as in the Brownian case using the $L^2$ isometry of compensated Poisson integrals. This is the main difficulty of the case $\beta <2$. In the one-dimensional case, Frikha and Li \cite{FrikhaWellposednessapproximationonedimensional2020} study a McKean-Vlasov
SDE driven by a compensated Poisson random measure with positive jumps. They give a rate of
convergence for the strong propagation of chaos in $L^1$ under one-sided Lipschitz assumptions on the coefficients. Let us mention that in dimension $d=1$, we obtain the same rates of convergence in \eqref{eq1_rq_exemple}. Concerning \eqref{eq3_rq_exemple} and \eqref{eq1_thm_poc_exemple}, our rate of convergence is of course worse than $N^{-\frac12}$ proved in the Brownian case in \cite{SznitmanTopicspropagationchaos1991}.

\section{Proof of  Theorem \ref{ThmexistenceMcKV1}}
Let us fix $\mu = (\mu_t)_{t\in [0,T]} \in \CC^0([0,T];\Pc_{\beta}(\R^d)).$  By using \cite[Theorem $6.2.9$]{ApplebaumLevyprocessesstochastic2009}, we deduce that the SDE \begin{equation}\label{existencemckvproof1}
	\left\{  \begin{array}{lll}
		&dX^{\mu}_t = b_t(X^{\mu}_t,\mu_t) \,dt + \sigma_t(X^{\mu}_{t^-},\mu_t)\, dZ_t, \quad t \in [0,T],\\ &X^{\mu}_0 = \xi \in L^{\beta}(\Omega,\FF_0;\R^d),
	\end{array}\right.
\end{equation}
admits a unique strong solution $X^{\mu}.$ Moreover, note that the coefficients of this standard SDE $(t,x) \in [0,T]\times \R^d \mapsto b_t(x,\mu_t)$ and $(t,x)\in [0,T]\times \R^d \mapsto \sigma_t(x,\mu_t)$ have at most linear growth with respect to the space variable $x,$ uniformly with respect to $t\in [0,T].$ By using Proposition $2$ in Fournier \cite{fournier_moment}, we get that $$\E \sup_{t \leq T} |X^{\mu}_t|^\beta <+\infty.$$ The map \begin{equation}
	\phi :\left\{  \begin{array}{rll}
		\CC^0([0,T];\Pc_{\beta}(\R^d)) & \rightarrow \CC^0([0,T];\Pc_{\beta}(\R^d)) \\ \mu &\mapsto ([X^{\mu}_t])_{t\in[0,T]}
	\end{array}\right.
\end{equation}
is thus well-defined. Indeed, the map $\phi$ takes values in $ \CC^0([0,T];\Pc_{\beta}(\R^d))$ by the dominated convergence theorem and by \cite[Lemma $2.3.2$]{ApplebaumLevyprocessesstochastic2009}, which ensures that for any $t \in [0,T]$, the process $X^{\mu}$ has no jump at time $t$ almost surely. The goal is now to prove that $\phi$ has a unique fixed point thanks to the Banach fixed point theorem. This is enough to prove the strong well-posedness of \eqref{Aedsmckv}. The space $\CC^0([0,T];\Pc_{\beta}(\R^d))$ is endowed with the uniform metric associated with $W_{\beta}.$ We fix $\mu,\nu \in  \CC^0([0,T];\Pc_{\beta}(\R^d))$ and we aim at estimating $\E \sup_{s\leq t}|X^{\mu}_s - X^{\nu}_s|^\beta,$ for $t \in [0,T].$ We employ the method used by Fournier in the proof of \cite[Proposition $2$]{fournier_moment}, which was used in the context of McKean-Vlasov SDEs by Frikha and Li in \cite{FrikhaWellposednessapproximationonedimensional2020} to prove the moment estimation \eqref{moment_estimate}. The first step is to consider the SDE \eqref{existencemckvproof1} without the big jumps term. Namely, we assume that for $\xi_1,\xi_2 \in L^\beta(\Omega,\FF_0;\R^d),$ and for all $t \in [0,T]$ \begin{equation}\label{existencemckvproof0}X^{\mu}_t = \xi_1 + \int_0^t b_s (X^{\mu}_s,\mu_s) \, ds + \int_0^t\int_{B_1} \sigma_s (X^{\mu}_{s^-}, \mu_s) z \, \NNN(ds,dz),\end{equation} 
and \begin{equation}\label{existencemckvproof00}X^{\nu}_t = \xi_2 + \int_0^t b_s (X^{\nu}_s,\nu_s) \, ds + \int_0^t\int_{B_1} \sigma_s (X^{\nu}_{s^-}, \nu_s) z \, \NNN(ds,dz).\end{equation}
Note that by definition of $\phi,$ $\xi_1$ is equal to $\xi_2,$ however in the next step of the proof, we need to take different initial data for the SDE. Using the Lipschitz assumption on the coefficients, Jensen and Burkholder-Davis-Gundy (BDG) inequalities, we obtain that for a constant $C = C_T$ depending only on $T$ and which can change from line to line, we have for all $t \in [0,T]$ \begin{multline*} \E \left( \sup_{s\leq t} |X^{\mu}_s - X^{\nu}_s|^2 \bigm|\xi_1,\xi_2 \right)\\ \leq  C \left[|\xi_1 - \xi_2|^2 + \int_0^t  \E \left( |X^{\mu}_s - X^{\nu}_s|^2 \bigm|\xi_1,\xi_2 \right)\, ds + \int_0^t W_{\beta}^2
	(\mu_s,\nu_s)\,ds\right].\end{multline*}  The conditional expectation with respect to the initial conditions in needed since they are not assumed to have a finite moment of order $2$. This is a technical point used in \cite{fournier:hal-00915365}. Grönwall's lemma ensures that $$ \E \left( \sup_{s\leq t} |X^{\mu}_s - X^{\nu}_s|^2 \bigm|\xi_1,\xi_2 \right)\leq  C|\xi_1 - \xi_2|^2  + C\int_0^t W_{\beta}^2
(\mu_s,\nu_s)\,ds.$$ It follows from Jensen's inequality that \begin{align*} \E \left( \sup_{s\leq t} |X^{\mu}_s - X^{\nu}_s|^\beta \bigm|\xi_1,\xi_2 \right)&\leq  \left(\E \left( \sup_{s\leq t} |X^{\mu}_s - X^{\nu}_s|^2 \bigm|\xi_1,\xi_2 \right) \right)^{\frac{\beta}{2}}\\ &\leq C|\xi_1 - \xi_2|^\beta  + C\left(\int_0^t W_{\beta}^2
	(\mu_s,\nu_s)\,ds\right)^{\frac{\beta}{2}}.\end{align*} Taking the expectation yields for all $t \in [0,T]$ \begin{equation}\label{existencemckvproof2}
	\E \left(\sup_{s\leq t} |X^{\mu}_s - X^{\nu}_s|^\beta\right) \leq C\E|\xi_1 - \xi_2|^\beta  + C\left(\int_0^t W_{\beta}^2
	(\mu_s,\nu_s)\,ds\right)^{\frac{\beta}{2}}.
\end{equation}
Let us now add the big jumps. We denote by $(T_n)_{n\geq1}$ the sequence of jumping times of $(Z_t)_{t\geq 0}$ having a size greater than $1,$ and by $(\Delta Z_n)_{n\geq 1}$ the associated sequence of jumps, which is an i.i.d.\ sequence of random variables with common distribution $\frac{\nu_{|B_1^c}}{\nu(B_1^c)}$ and independent of $(T_n)_{n\geq 1}.$ We can write the restriction of the Poisson random measure $\NN$ on $\R^+ \times B_1^c$ as $$ \sum_{n \geq 1} \delta_{(T_n,\Delta Z_n)},$$ which is independent of the restriction of $\NN$ on $\R^+ \times B_1\backslash\{0\}.$ Let us denote by $\mathcal{G}$ the $\sigma$-algebra generated by $(T_n)_{n\geq 1}.$ Notice that on the time interval $[0,T_1),$ $X^{\mu}$ and $X^{\nu}$ defined in \eqref{existencemckvproof1} are respectively solutions to \eqref{existencemckvproof0} and \eqref{existencemckvproof00} with $\xi_1 = \xi_2.$ Thus, using \eqref{existencemckvproof2} with the conditional expectation with respect to $\mathcal{G}$ instead of the expectation, we deduce that \begin{equation}\label{existencemckvproof3}
	\E \left(\sup_{s < t \wedge T_1} |X^{\mu}_s - X^{\nu}_s| ^\beta \bigm| \mathcal{G} \right)\leq C \left(\int_0^{t\wedge T_1} W_{\beta}^{2}(\mu_s,\nu_s)\, ds\right)^{\frac{\beta}{2}}.
\end{equation} Let us now deal with the first big jump of $Z$, which occurs at time $T_1.$ We have  $$X^{\mu}_{T_1} - X^{\nu}_{T_1} = X^{\mu}_{T_1^-} - X^{\nu}_{T_1^-}  + \left(\sigma_{T_1}(X^{\mu}_{T_1^-},\mu_{T_1}) -\sigma_{T_1}(X^{\nu}_{T_1^-},\nu_{T_1})\right)  \Delta Z_1.$$ It follows from the Lipschitz assumption on $\sigma$ that $$ | X^{\mu}_{T_1} - X^{\nu}_{T_1} |^\beta \leq C|X^{\mu}_{T_1^-} - X^{\nu}_{T_1^-}|^\beta ( 1 + |\Delta Z_1|^\beta) + W_{\beta}^{\beta}(\mu_{T_1},\nu_{T_1})|\Delta Z_1|^\beta.$$ Since $\Delta Z_1$ is independent of $\mathcal{G}$ and $\E|\Delta Z_1|^\beta <+ \infty,$ we deduce by \eqref{existencemckvproof3} that conditionally on the event $\{T_1 \leq t \}$ $$ \E \left(| X^{\mu}_{T_1} - X^{\nu}_{T_1} |^\beta \bigm| \mathcal{G}\right) \1_{ T_1 \leq t} \leq C \left[ \left(\int_0^{t\wedge T_1} W_{\beta}^{2}(\mu_s,\nu_s) \, ds\right)^{\frac{\beta}{2}} + W_{\beta}^\beta (\mu_{T_1}, \nu_{T_1}) \right].$$ We thus have by the preceding inequality and \eqref{existencemckvproof3} \begin{align*}
	&\E \left(| X^{\mu}_{t\wedge T_1} - X^{\nu}_{t\wedge T_1} |^\beta \bigm| \mathcal{G}\right)\\&= 	\E \left(| X^{\mu}_{T_1} - X^{\nu}_{ T_1} |^\beta \bigm| \mathcal{G}\right) \1_{T_1 \leq t} + \E \left( | X^{\mu}_{t} - X^{\nu}_{t} |^\beta \bigm| \mathcal{G}\right) \1_{T_1 >t}  \\ &\leq C \left[ \left(\int_0^{t\wedge T_1} W_{\beta}^{2}(\mu_s,\nu_s) \, ds\right)^{\frac{\beta}{2}} + W_{\beta}^\beta (\mu_{T_1}, \nu_{T_1}) \right].
\end{align*}

Following the same lines and using \eqref{existencemckvproof2}, we prove that for any $n \geq 1$ \begin{multline}\label{existencemckvproof4}
	\E \left(\sup_{ t \wedge T_n \leq s < t \wedge T_{n+1}}| X^{\mu}_{s} - X^{\nu}_{s} |^\beta \bigm| \mathcal{G}\right)\\ \leq 	C\left[\E \left(| X^{\mu}_{t\wedge T_n} - X^{\nu}_{t\wedge T_n} |^\beta \bigm| \mathcal{G}\right) +  \left(\int_{t\wedge T_n}^{t \wedge T_{n+1}} W_{\beta}^2(\mu_s,\nu_s) \, ds \right)^{\frac{\beta}{2}}\right],
\end{multline}
and \begin{align}\label{existencemckvproof5}
	&\E \left(\sup_{ t \wedge T_n \leq s \leq t \wedge T_{n+1}}| X^{\mu}_{s} - X^{\nu}_{s} |^\beta \bigm| \mathcal{G}\right)\\ \notag & \hspace{0.5cm} \leq 	C \left[\E \left(| X^{\mu}_{t\wedge T_n} - X^{\nu}_{t\wedge T_n} |^\beta \bigm| \mathcal{G}\right) +  \left(\int_{t\wedge T_n}^{t \wedge T_{n+1}} W_{\beta}^2(\mu_s,\nu_s) \, ds \right)^{\frac{\beta}{2}} + W_{\beta}^2(\mu_{T_{n+1}}, \nu_{T_{n+1}}) \right],
\end{align}
where $C$ is independent of $n$. Reasoning by induction, we deduce that for a constant $C>1$ depending only on $T,$ we have \begin{align*}
	\notag &\E \left(\sup_{ t \wedge T_n \leq s < t \wedge T_{n+1}}| X^{\mu}_{s} - X^{\nu}_{s} |^\beta \bigm| \mathcal{G}\right) \\&\leq 	C^{n+1}\left[  \sum_{k=0}^n\left(\int_{t\wedge T_k}^{t \wedge T_{k+1}} W_{\beta}^2(\mu_s,\nu_s) \, ds \right)^{\frac{\beta}{2}} + \sum _{k=1}^n W_{\beta}^\beta(\mu_{T_k},\nu_{T_k})\right] \\ &\leq C^{n+1}\left[  (n+1)^{1 - \frac{\beta}{2}}\left(\int_{0}^{t \wedge T_{n+1}} W_{\beta}^2(\mu_s,\nu_s) \, ds \right)^{\frac{\beta}{2}} + \sum _{k=1}^n W_{\beta}^\beta(\mu_{T_k},\nu_{T_k})\right],
\end{align*}
by Jensen's inequality and with the convention that $T_0 = 0.$ Thus, for a certain constant $ K >C,$ one has for all $ j \leq n$  \begin{multline*}
	\notag \E \left(\sup_{ t \wedge T_j \leq s < t \wedge T_{j+1}}| X^{\mu}_{s} - X^{\nu}_{s} |^\beta \bigm| \mathcal{G}\right) \\\leq K^{j+2}\left[ \left(\int_{0}^{t \wedge T_{n+1}} W_{\beta}^2(\mu_s,\nu_s) \, ds \right)^{\frac{\beta}{2}} + \sum _{k=1}^n W_{\beta}^\beta(\mu_{T_k},\nu_{T_k})\right].
\end{multline*}Summing the preceding inequality over $j \in \{0,\dots,n\}$, we deduce that  \begin{multline}\label{existencemckvproof7}
	\E \left(\sup_{ 0\leq s < t \wedge T_{n+1}}| X^{\mu}_{s} - X^{\nu}_{s} |^\beta \bigm| \mathcal{G}\right) \\ \leq \frac{K^{n+3}}{1-K}\left[ \left(\int_{0}^{t \wedge T_{n+1}} W_{\beta}^2(\mu_s,\nu_s) \, ds \right)^{\frac{\beta}{2}} + \sum _{k=1}^n W_{\beta}^\beta(\mu_{T_k},\nu_{T_k})\right].
\end{multline} 
Let us denote by $(N_t)_{t \geq 0}$ the Poisson process associated with the jumping times $(T_n)_{n\geq 1}$ which has an intensity $\lambda = \nu(B_1^c).$  One has \begin{align*}
	&\E \left(\sup_{ 0\leq s \leq t}| X^{\mu}_{s} - X^{\nu}_{s} |^\beta \right) \\&= \E \left( \E \left(\sup_{ 0\leq s \leq t }| X^{\mu}_{s} - X^{\nu}_{s} |^\beta \bigm| \mathcal{G}\right) \1_{ t < T_{1}}\right)  + \sum_{n=1}^{\infty} \E \left( \E \left(\sup_{ 0\leq s \leq t }| X^{\mu}_{s} - X^{\nu}_{s} |^\beta \bigm| \mathcal{G}\right) \1_{T_n \leq t < T_{n+1}}\right) \\ &\leq \E \left( \E \left(\sup_{ 0\leq s < t \wedge T_1 }| X^{\mu}_{s} - X^{\nu}_{s} |^\beta \bigm| \mathcal{G}\right) \1_{ t < T_{1}}\right) \\ &\quad + \sum_{n=1}^{\infty} \E \left( \E \left(\sup_{ 0\leq s < t \wedge T_{n+1} }| X^{\mu}_{s} - X^{\nu}_{s} |^\beta \bigm| \mathcal{G}\right) \1_{T_n \leq t < T_{n+1}}\right)
\end{align*}

Using \eqref{existencemckvproof3} and \eqref{existencemckvproof7}, we obtain that for any $t \in [0,T]$ \begin{align}\label{existencemckvproof8}
	&\notag\E \left(\sup_{ 0\leq s \leq t}| X^{\mu}_{s} - X^{\nu}_{s} |^\beta \right)\\ \notag&\leq C \left(\int_0^{t} W_{\beta}^{2}(\mu_s,\nu_s)\right)^{\frac{\beta}{2}} \\ \notag &\quad + \sum_{n=1}^{\infty} \frac{K^{n+3}}{1-K} \E \left(\left[ \left(\int_{0}^{t \wedge T_{n+1}} W_{\beta}^2(\mu_s,\nu_s) \, ds \right)^{\frac{\beta}{2}} + \sum _{k=1}^n W_{\beta}^\beta(\mu_{T_k},\nu_{T_k})\right] \1_{T_n \leq t < T_{n+1}}\right) \\ &\leq   C \left(\int_0^{t} W_{\beta}^{2}(\mu_s,\nu_s)\right)^{\frac{\beta}{2}} \\ \notag &\quad+ \sum_{n=1}^{\infty} \frac{K^{n+3}}{1-K} \P(N_t =n) \left( \left(\int_{0}^{t} W_{\beta}^2(\mu_s,\nu_s) \, ds \right)^{\frac{\beta}{2}} + \E \left(\sum _{k=1}^n W_{\beta}^\beta(\mu_{T_k},\nu_{T_k}) \bigm| N_t=n \right)\right).
\end{align}
Let us recall that the conditional distribution of $(T_1, \dots, T_n)$ given $N_t=n$ admits the following density with respect to the Lebesgue measure $$ (t_1, \dots, t_n) \in [0,t]^n \mapsto \frac{n!}{t^n} \1_{ t_1 < \dots < t_n}.$$ This yields \begin{align*}
	\E \left(\sum _{k=1}^n W_{\beta}^\beta(\mu_{T_k},\nu_{T_k}) \bigm| N_t=n \right) &= \int_{[0,t]^n}  \frac{n!}{t^n} \1_{t_1 < \dots < t_n} \sum _{k=1}^n W_{\beta}^\beta(\mu_{t_k},\nu_{t_k}) \, dt_1 \dots\, dt_n \\\notag &= \int_{[0,t]^n}  \frac{1}{t^n} \sum _{k=1}^n W_{\beta}^\beta(\mu_{t_k},\nu_{t_k})\, dt_1 \dots\, dt_n \\\notag &= \sum_{k=1}^n \frac{1}{t^n} t^{n-1}\int_0^t W_{\beta}^\beta(\mu_s,\nu_s) \, ds \\\notag &= \frac{n}{t} \int_0^t W_{\beta}^\beta(\mu_s,\nu_s) \, ds.
\end{align*}
Injecting this equality in \eqref{existencemckvproof8}, we get \begin{align*}
	&\E \left(\sup_{ 0\leq s \leq t}| X^{\mu}_{s} - X^{\nu}_{s} |^\beta \right)\\ \notag&\leq   \sum_{n=1}^{\infty} \frac{K^{n+3}}{1-K} \frac{(\lambda t)^n}{n!}e^{-\lambda t} \left( \left(\int_{0}^{t} W_{\beta}^2(\mu_s,\nu_s) \, ds \right)^{\frac{\beta}{2}} +\frac{n}{t} \int_0^t W_{\beta}^\beta(\mu_s,\nu_s) \, ds\right) \\ &\quad+ C \left(\int_0^{t} W_{\beta}^{2}(\mu_s,\nu_s)\right)^{\frac{\beta}{2}} .
\end{align*}
This proves the existence of a constant $C>0$ depending only on $T$ such that for all $t \in [0,T]$

\begin{equation}\label{existencemckvproof9bis}
	\E \left(\sup_{ 0\leq s \leq t}| X^{\mu}_{s} - X^{\nu}_{s} |^\beta \right) \leq  C \left[\left(\int_0^{t} W_{\beta}^{2}(\mu_s,\nu_s) \, ds\right)^{\frac{\beta}{2}} + \int_0^t W_{\beta}^\beta(\mu_s,\nu_s) \, ds\right].
\end{equation}
 Changing again the constant $C,$ Hölder's inequality yields for all $t \in [0,T]$ 
\begin{equation}\label{existencemckvproof9}
	\sup_{0\leq s \leq t}W_{\beta}^\beta(\phi(\mu)_s,\phi(\nu)_s) \leq\E \left(\sup_{ 0\leq s \leq t}| X^{\mu}_{s} - X^{\nu}_{s} |^\beta \right) \leq  C \left(\int_0^{t} W_{\beta}^{2}(\mu_s,\nu_s) \, ds\right)^{\frac{\beta}{2}}.
\end{equation}
Raising the preceding inequality to the power $\frac{2}{\beta}$ and reasoning by induction, we prove that for any $ n\geq 1$ and for any $t \in [0,T]$
\begin{equation*}
	\sup_{0\leq s \leq t}W_{\beta}^\beta(\phi^n(\mu)_s,\phi^n(\nu)_s) \leq C^n\left(\frac{t^n}{n!}\right)^{\frac{\beta}{2}}\sup_{0 \leq s \leq t} W_{\beta}^\beta (\mu_s, \nu_s).
\end{equation*}
This proves that for $n$ large enough, $\phi^n$ is a contraction on $\CC^0([0,T];\Pc_{\beta}(\R^d))$. The function $\phi$ has thus a unique fixed point by the Banach fixed point theorem, which concludes the proof.
%

\section{Proof of Theorem \ref{thm_poc_lip}}

\noindent\textbf{Proof of \eqref{eq1_rq_exemple} and \eqref{eq2_rq_exemple}.} We write for all $t \in [0,T]$ \begin{align*}  X^{i,N}_{t} - X^{i,\infty}_{t}&= \int_0^t b_s(X^{i,N}_{s},\muu^N_{s}) - b_s(X^{i,\infty}_{s},\mu_{s}) \, ds.
\end{align*}
Using the Lipschitz assumption on $b$, there exists $C>0$ such that for all $t \in [0,T]$ \begin{align*}
	&\sup_{i \leq N}\E  \sup_{r \leq t} \vert X^{i,N}_{r} - X^{i,\infty}_{r} \vert\\ &\leq C \int_0^t \sup_{i \leq N}\E \vert X^{i,N}_{s} - X^{i,\infty}_{s} \vert \, ds + C\int_0^t \E W_1(\muu^N_{s},\mu_{s}) \, ds \\ &\leq C \int_0^t \sup_{i \leq N}\E  \vert X^{i,N}_{s} - X^{i,\infty}_{s} \vert \, ds + C\int_0^t \E W_1(\muu^N_{s},\tilde{\mu}^N_{s}) + \E W_1(\tilde{\mu}^N_{s},\mu_{s}) \, ds,
\end{align*}
where $\tilde{\mu}^N_{s}:= \frac{1}{N} \sum\limits_{k=1}^N \delta_{X^{k,\infty}_{s}}$ is the empirical measure associated with $(X^{i,\infty})_{i\geq 1}$. Using that $$W_1(\muu^{N}_{s},\tilde{\mu}^N_{s}) \leq \frac{1}{N}\sum_{k=1}^N|X^{k,N}_{s} - X^{k,\infty}_{s}|$$ and Grönwall's inequality, there exists $C_T>0$, depending on the final horizon of time $T$, such that for all $N\geq 1$ \begin{equation}\label{control_Gronwall_lip}
	\sup_{i \leq N}	\E  \sup_{t \leq T} \vert X^{i,N}_{t} - X^{i,\infty}_{t} \vert \leq C \int_0^T \E W_1(\tilde{\mu}^N_{s},\mu_{s}) \, ds.\end{equation}
We conclude using \cite[Theorem $1$]{fournier:hal-00915365} since $(X^{i,\infty})_{i \geq 1}$ are i.i.d.\ and $$  \sup_{i \geq 1} \sup_{t \in[0,T]}\E |X^{i,\infty}_{t}|^\beta <+ \infty$$ by Grönwall's inequality. The inequality \eqref{eq2_rq_exemple} follows from \eqref{eq1_rq_exemple} and \cite{fournier:hal-00915365} because \begin{align*}
	\sup_{t \in[0,T]} \E W_1(\muu^N_t,\mu_t) & \leq  \sup_{t \in[0,T]} \E W_1(\muu^N_t,\tilde{\mu}^N_t) + \sup_{t \in[0,T]} \E W_1(\tilde{\mu}^N_t,\mu_t)\\ &\leq \sup_{t\in[0,T]}  \sup_{i\leq N}\E |X^{i,N}_t - X^{i,\infty}_t| + \sup_{t \in[0,T]} \E W_1(\tilde{\mu}^N_t,\mu_t).
\end{align*}

\noindent\textbf{Proof of \eqref{eq3_rq_exemple}.} As in the proof of \cite[Theorem $1.4$]{SznitmanTopicspropagationchaos1991}, we  write for all $t \in [0,T]$ \begin{align*}  X^{i,N}_{t} - X^{i,\infty}_{t}&= \int_0^t b_s(X^{i,N}_{s},\muu^N_{s}) - b_s(X^{i,\infty}_{s},\mu_{s}) \, ds \\ &= \int_0^t \frac1N \sum_{j=1}^N \Big\{ \tilde{b}(X^{i,N}_{s},X^{j,N}_s) - \tilde{b}(X^{i,\infty}_{s},X^{j,N}_s) + \tilde{b}(X^{i,\infty}_{s},X^{j,N}_s) - \tilde{b}(X^{i,\infty}_{s},X^{j,\infty}_s) \\ & \hspace{2cm} +  \tilde{b}(X^{i,\infty}_{s},X^{j,\infty}_s) - \int_{\R^d} \tilde{b}(X^{i,\infty}_s,y) \, d\mu_s(y)  \Big\} \, ds.
\end{align*}
Let us set, for $t \in [0,T]$ and $x,x' \in \R^d$, $$B_t(x,x') :=  \tilde{b}(x,x') - \int_{\R^d} \tilde{b}(x,y) \, d\mu_t(y).$$
Reasoning exactly as in the proof of \cite[Theorem $1.4$]{SznitmanTopicspropagationchaos1991}, i.e.\ using the exchangeability of the particles, the Lipschitz continuity of $\tilde{b}$ and Grönwall's inequality, we obtain that \begin{align*}
	&\E \sup_{t \in [0,T]} |X^{i,N}_t - X^{i, \infty}_t| \\ &\leq C \int_0^T \E \Big(\frac1N \sum_{j=1}^N |B_s(X^{i,\infty}_s,X^{j,\infty}_s)|  \Big) \, ds\\ &\leq  \frac{C}{N} \left[
	\int_0^T  \E|B_s(X^{i,\infty}_s,X^{i,\infty}_s)+ \int_{\R^d} \Big(\E \Big(\sum_{j\neq i} |B_s(x,X^{j,\infty}_s)| \Big)^\beta \Big)^{\frac1\beta} \, d\mu_s(x) \, ds\right].
\end{align*}
By the $\beta$-th order Marcinkiewicz-Zygmund inequality (see $(1.4)$ in \cite{ChenSunginequalitysum}), we deduce that 
\begin{multline*}
	\E \sup_{t \in [0,T]} |X^{i,N}_t - X^{i, \infty}_t|  \\
	\leq   \frac{C}{N}\left[ \int_0^T \E|B_s(X^{i,\infty}_s,X^{i,\infty}_s)|  +\int_{\R^d} \Big(\sum_{j\neq i} \E |B_s(x,X^{j,\infty}_s)|^\beta\Big)^{\frac1\beta}  d\mu_s(x) \, ds\right].
\end{multline*}
Using that $B_s$ has at most linear growth on $\R^d \times \R^d$, uniformly with respect to $s \in [0,T]$ by \eqref{moment_estimate}, we obtain that 

\begin{align}\label{eq_conclusion}
	\E \sup_{t \in [0,T]} |X^{i,N}_t - X^{i, \infty}_t|  & \leq C N^{\frac1\beta - 1} (\E \sup_{t\in[0,T]} |X^{1,\infty}_t|^\beta)^{\frac{1}{\beta}}  \\\notag  &\leq CN^{\frac1\beta - 1},
\end{align}
where the last inequality follows from \eqref{moment_estimate}.\\

\noindent \textbf{Proof of \eqref{eq1_thm_poc_exemple}.} Let us recall that since $Z^i$ is a $\alpha$-stable process with $\alpha \in (1,2)$, it can be written as $$Z^i_t = \int_0^t\int_{\R^d} z \NNN^i(ds,\,dz),$$ see Remark $14.6$ and Theorem $14.7$ in \cite{SatoLevyprocessesinfinitely1999}. As a first step, we remove the jumps of size larger than the number of particles $N$ from all the noises, i.e.\ we define, for $i \geq 1$ and $t \in [0,T]$ $$ Z^i_{N,t} := \int_0^t \int_{B_N} z \, \NNN^i(ds,dz).$$ We define, for all $i \in \{1,\dots,N\}$, $X^{i,\infty}_N$ as the unique solution to

\begin{equation}\label{edsmckvtronquee1_exemple}
	\begin{cases}
		&dX^{i,\infty}_{N,t}  = b_t(X^{i,\infty}_{N,t}, \mu_{N,t})\, dt +\,dZ^i_{N,t}, \quad t \in [0,T],\quad i \in \{1,\dots,N \}, \\ &\mu_{N,t} := [X^{i,\infty}_{N,t}],\\ &X^{i,\infty}_{N,0}  = \xi^i.
	\end{cases}
\end{equation}
For any $N \geq 1$ fixed, the random variables $(X^{i,\infty}_{N})_{ i \leq N}$ are i.i.d. We proceed similarly for the particle system by defining $(X^{i,N}_{N})_{i\leq N}$ as the unique solution to  \begin{equation}\label{edsparticles_tronquee_exemple}
	\begin{cases}
		&dX^{i,N}_{N,t} = b_s(X^{i,N}_{N,t}, \muu^N_{N,t}) \,dt+ \,dZ^i_{N,t}, \quad t \in [0,T],\quad i \in \{1,\dots,N\}, \\ &\muu^N_{N,t} := \frac{1}{N} \sum\limits_{j=1}^N \delta_{X^{j,N}_{N,t}},\\ &X^{i,N}_{N,0} = \xi^i,
	\end{cases}\end{equation} 
The first objective is to control the  $L^1$-error respectively between $X^{i,N}_{N}$ and $X^{i,N}$ and between $X^{i,\infty}_N$ and $X^{i,\infty}$ for all $i \in \{1,\dots,N\}.$ We write for all $t \in [0,T]$ $$ X^{i,N}_{N,t} - X^{i,N}_t = \int_0^t (b_s(X^{i,N}_{N,s}, \muu^N_{N,s}) - b_s(X^{i,N}_s,\muu^N_s))\, ds - \int_0^t \int_{B_N^c} z \, \NNN^i(ds,dz).$$  Using the fact that $b$ is Lipschitz continuous on $\R^d \times \mathcal{P}_1(\R^d)$ and BDG's inequality, there exists $C>0$ independent of $N\geq 1$ and $t \in[0,T]$, which can change from line to line, such that for all $t \in [0,T]$\begin{multline*}
	\sup_{i \leq N}\E  \sup_{r\leq t} |X^{i,N}_{N,r} - X^{i,N}_r |
	\leq C \left[ \int_0^t\sup_{i \leq N}\E|X^{i,N}_{N,s} - X^{i,N}_s |\, ds + \int_0^t \frac1N \sum_{j=1}^N \E|X^{j,N}_{N,s} - X^{j,N}_s |\, ds \right.\\ \left.   + \sup_{i\leq N}\E \left(\int_0^t\int_{B_N^c} |z|^2 \, \NN^i(ds,dz)\right)^{\frac12}\right].
\end{multline*} Using the subadditivity of the square root, one has for all $t \in [0,T]$ 	\begin{align*}\sup_{i \leq N}\E  \sup_{r\leq t} |X^{i,N}_{N,r} - X^{i,N}_r |
	&\leq C \left[ \int_0^t\sup_{i \leq N}\E|X^{i,N}_{N,s} - X^{i,N}_s |\, ds + \int_0^t\int_{B_N^c} |z| \, d\nu(z)\,ds\right] \\ &\leq C \left[ \int_0^t\sup_{i \leq N}\E|X^{i,N}_{N,s} - X^{i,N}_s |\, ds + \int_{N}^{\infty} r \, \frac{dr}{r^{1+\alpha}}\right].
\end{align*} Grönwall's inequality ensures that \begin{equation}\label{ThmPC3proofeqi_exemple}
	\sup_{i \leq N}\E \sup_{t\leq T} |X^{i,N}_{N,t} - X^{i,N}_t | \leq CN^{1-\alpha}.
\end{equation} We similarly get that for some constant $C>0$ independent of $N$
\begin{equation}\label{ThmPC3proofeqiii_exemple}
	\sup_{i \leq N}\E \sup_{t\leq T}|X^{i,\infty}_{N,t} - X^{i,\infty}_t| \leq CN^{1-\alpha}.
\end{equation}
The triangle inequality, \eqref{ThmPC3proofeqi_exemple} and \eqref{ThmPC3proofeqiii_exemple} yield \begin{equation}\label{erreur_eq_tronquee_exemple}	\sup_{i \leq N}\E \sup_{t \leq T} |X^{i,N}_t - X^{i,\infty}_t| \leq CN^{1-\alpha} + 	\sup_{i \leq N}\E  \sup_{t \leq T} |X^{i,N}_{N,t} - X^{i,\infty}_{N,t}|.\end{equation}

Let us control the second term in the right hand-side of \eqref{erreur_eq_tronquee_exemple}. Reasoning as in the proof of \ref{thm_poc_lip} (see \eqref{eq_conclusion}), we have, for $\beta = \alpha$, 

\begin{align}\label{eq_conclusion_2}
	\E \sup_{t \in [0,T]} |X^{i,N}_{N,t} - X^{i, \infty}_{N,t}|  & \leq C N^{\frac1\alpha - 1} (\E \sup_{t\in[0,T]} |X^{1,\infty}_{N,t}|^\alpha)^{\frac{1}{\alpha}}.
\end{align}
Grönwall's inequality ensures that there exists $C>0$ such that for any $t \in [0,T]$  $$\E \sup_{s\leq t}|X^{i,\infty}_{N,s}|^\alpha \leq  C \sup_{s \leq t}\E |Z^{i}_{N,s}|^\alpha.$$ Then, BDG's and Jensen's inequalities yield \begin{align}\label{moment_sup_alpha_exemple}
	\notag\E|Z^{i}_{N,t}|^\alpha &\leq C \E \left(\int_0^t \int_{B_N} |z|^2 \, \NN^i(ds,dz) \right)^{\frac{\alpha}{2}}\\ \notag &\leq C \left[\left(\E \int_0^t \int_{ B_1} |z|^2 \,\NN^i(ds,dz) \right)^{\frac{\alpha}{2}}  + \E \left(\int_0^t \int_{B_N\backslash B_1} |z|^2 \, \NN^i(ds,dz) \right)^{\frac{\alpha}{2}} \right]\\ &\leq C \left[1+ \int_1^N r^\alpha \,\frac{dr}{r^{1+\alpha}} \right]\\ \notag&\leq C \ln(N).
\end{align} 
Coming back to \eqref{eq_conclusion_2}, we have \begin{align*}
	\E \sup_{t \in [0,T]} |X^{i,N}_{N,t}- X^{i, \infty}_{N,t}|  & \leq C N^{\frac1\alpha - 1} (\ln(N))^{\frac{1}{\alpha}}.
\end{align*}
Since $\alpha >1$ and thus $\frac1\alpha - 1 < 1 - \alpha<0$, we deduce from \eqref{erreur_eq_tronquee_exemple} that $$ \E \sup_{t \in [0,T]}|X^{i,N}_t - X^{i,\infty}_t| \leq C N^{\frac1\alpha - 1} (\ln(N))^{\frac{1}{\alpha}},$$ which ends the proof.

\appendix
	
	\section{Existence of a solution to \eqref{Aedsmckv} when $\beta \in (0,1)$}\label{appendix_existence}
	
	Let us fix $\beta \in (0,1)$. We have seen in Remark \ref{_rq_beta_petit} that in this case, uniqueness for \eqref{Aedsmckv} fails to be true under Assumption \ref{Assumpt_lip}. However, the existence of solutions to \eqref{Aedsmckv} is given in the following proposition.
	
	\begin{Prop}\label{ThmexistenceMcKV1bis}
		We assume that Assumption \ref{Assumpt_lip} is satisfied and that there exists $\delta >0$ such that $$ \int_{B_1^c} |z|^{\beta + \delta} \, d\nu(z) <+ \infty.$$ Then, there exists a strong solution $(X_t)_{t \in [0,T]}$ to \eqref{Aedsmckv} for all $\xi \in L^{\beta + \delta}(\Omega,\FF_0;\R^d).$ Moreover, we have \begin{equation}\label{thmexistencemckvbismoment}
			\E \sup_{t\leq T} |X_t|^{\beta+ \delta}  <+ \infty.	\end{equation}
	\end{Prop}
	
	\begin{proof}
		The strategy relies on a compactness argument. Let us denote by $D_T:=D([0,T];\R^d)$ the Skorokhod space, i.e. the space of càdlàg $\R^d$-valued functions defined on $[0,T].$ We endow $D_T$ with the Skorokhod metric $d,$ which makes it Polish (see \cite[Section $34$]{Bass_stochastic}). By definition of $d,$ we have for any $f,g \in D_T$ $$ d(f,g) \leq \Vert f-g \Vert_{\infty} := \sup_{t \in [0,T]} |f_t - g_t|.$$ The previous inequality becomes an equality if $g=0.$ We also denote by $\Pc_{\beta}(D_T)$ the space of probability measures $\mu \in \Pc(D_T)$ such that $$ \int_{D_T} d(f,0)^{\beta} \, d\mu(f) =\int_{D_T} \Vert f\Vert_{\infty}^{\beta} \, d\mu(f) < + \infty.$$ It is endowed with the Wasserstein metric of order $\beta$ defined, for any $\mu,\nu \in \Pc_{\beta}(D_T)$,  by $$\mathcal{W}_{\beta}(\mu,\nu) = \inf_{\pi \in \Pi(\mu,\nu)} \int_{D_T\times D_T} d^{\beta}(f,g) \, d\pi(f,g),$$ where $\Pi(\mu,\nu)$ denotes the set of probability measure on $D_T \times D_T$ having $\mu$ and $\nu$ as marginal distributions. For any fixed $t \in [0,T],$ we define the projection $\pi_t : f\in D_T \mapsto f_t \in \R^d.$ It is a measurable function so that if $\mu$ belongs to $\Pc(D_T),$ we can define $\mu_t \in \mathcal{P}(\R^d)$ as the push-forward measure of $\mu$ by $\pi_t.$ Notice that if $\mu \in \Pc_{\beta}(D_T),$ the function $t \in [0,T] \mapsto \mu_t $ belongs to $D([0,T];\Pc_{\beta}(\R^d)).$ Let us fix $\mu \in \Pc_{\beta}(D_T).$ Reasoning as in the proof of Proposition \ref{ThmexistenceMcKV1}, the standard SDE \begin{equation}\label{existencemckvbisproof1}
			\left\{  \begin{array}{lll}
				&dX^{\mu}_t = b_t(X^{\mu}_t,\mu_t) \,dt + \sigma_t(X^{\mu}_{t^-},\mu_t)\, dZ_t, \quad t\in[0,T], \\ &X^{\mu}_0 = \xi \in L^{\beta}(\Omega,\FF_0)
			\end{array}\right.
		\end{equation}
		admits a unique strong solution $X^{\mu}$  such that $$\E \sup_{t \leq T} |X^{\mu}_t|^\beta <+\infty.$$ The following function is thus well-defined \begin{equation}
			\phi :\left\{  \begin{array}{rll}
				\Pc_{\beta}(D_T) & \rightarrow\Pc_{\beta}(D_T) \\ \mu &\mapsto [(X^{\mu}_t)_{t\in[0,T]}].
			\end{array}\right.
		\end{equation}
		The goal is to prove that $\phi$ has at least one fixed point using Schauder's fixed point theorem. By the estimation \eqref{existencemckvproof9bis} obtained in the proof of Theorem \ref{ThmexistenceMcKV1}, we have
		
		\begin{equation}\label{existencemckvbisproof2}
			\mathcal{W}_{\beta}(\phi(\mu),\phi(\nu))\leq\E \left(\sup_{ 0\leq s \leq T}| X^{\mu}_{s} - X^{\nu}_{s} |^\beta \right) \leq  C \left(\int_0^{T} W_{\beta}^{2}(\mu_s,\nu_s) \, ds\right)^{\frac{\beta}{2}}.
		\end{equation}
		Let us show that this implies the continuity of $\phi.$ Consider $(\mu^n)_n \in \Pc_{\beta}(D_T)$ a sequence which converges towards $\mu \in \Pc_{\beta}(D_T)$ with respect to $\mathcal{W}_{\beta}.$ For almost all $t \in [0,T],$ the sequence $(\mu^n_t)_n$ converges in distribution to $\mu_t$ (see \cite[Section $13$]{BillingsleyConvergenceprobabilitymeasures1999}). Let us fix such a $t \in [0,T].$ We aim at proving that the previous convergence holds true with respect to $W_{\beta}.$ By \cite[Definition $6.8$]{Villaniold}, it is enough to prove that $$ \lim_{R \rightarrow + \infty} \sup_{n} \int_{|x|\geq R} |x|^{\beta} d\mu^n_t(x)= 0.$$
		But since \begin{align*}
			\int_{|x|\geq R} |x|^{\beta} d\mu^n_t(x) &= \int_{|f_t|\geq R} |f_t|^{\beta} d\mu^n(f) \\ &\leq  \int_{d(f,0)\geq R} d(f,0)^{\beta} d\mu^n(f),
		\end{align*}
		we conclude using that $\mu^n \overset{\mathcal{W}_{\beta}}{\longrightarrow} \mu.$
		Thus, for almost all $t \in [0,T],$ $(\mu^n_t)_n$ converges to $\mu_t$ with respect to $W_{\beta}.$ Coming back to \eqref{existencemckvbisproof2}, and using the dominated convergence theorem justified since $$ \sup_{n\geq 1} \int_{D_T} \Vert f \Vert_{\infty}^\beta \, d\mu^n(f) < + \infty,$$ we conclude that $\phi$ is continuous. Following the same lines as to prove \eqref{existencemckvproof9bis}, we show that for some constant $C=C_T>0$ \begin{equation}\label{existencemckvbisproof3}
			\E \sup_{t \leq T} |X^{\mu}_t|^\beta \leq C \left[1+\left(\int_0^T M_{\beta}^2(\mu_s) \, ds  \right)^{\frac{\beta}{2}}\right].	\end{equation} Let us define, for $R>0$,  $$ \overline{\mathcal{B}}_R := \left\{ \mu \in \Pc_{\beta}(D_T), \, \int_{D_T} \Vert f \Vert_{\infty}^\beta \, d\mu(f) \leq R\right\}.$$ This is a closed and convex subset of $\Pc_{\beta}(D_T),$ which is stable by $\phi$ for $R$ large enough owing to \eqref{existencemckvbisproof3} and since $\beta \in (0,1)$. In the following, we fix $R>0$ such that $\phi(\overline{\mathcal{B}}_R) \subset \overline{\mathcal{B}}_R.$ It remains to prove that $\phi(\overline{\mathcal{B}}_R)$ is relatively compact in $\Pc_{\beta}(D_T)$ to conclude that $\phi$ admits a fixed point by Schauder's theorem. Let us fix $(\mu^n)_n$ a sequence of $\overline{\mathcal{B}}_R.$ In a first step, we prove with Aldou's criterion (see \cite[Theorem $34.8$]{Bass_stochastic}) that $([(X^{\mu^n}_t)_{t\in[0,T]}])_n$ is tight, and thus relatively compact in $\Pc(D_T).$ For $t \in [0,T]$ and $A>0,$ we have \begin{align*}
			\P(|X^{\mu^n}_t| \geq  A) &\leq \frac{1}{A^{\beta}} \E |X^{\mu^n}_t|^{\beta} \\ &\leq \frac{R}{A^{\beta}}. \end{align*}
		This yields  for all $t \in [0,T]$ \begin{equation}\label{Aldous1}
			\lim_{A \rightarrow + \infty} \sup_{n}\, \P(|X^{\mu^n}_t| \geq A) = 0.
		\end{equation}
		Let $(\tau_n)_n$ be a sequence of stopping times and $(\delta_n)_n $ a sequence of real numbers converging to $0.$ We assume that $\tau_n \leq T$ and $0 \leq \tau_n + \delta_n \leq T$ for any $n \geq 1$. It remains to prove that for any fixed $\epsilon >0$ $$ \P (|X^{\mu^n}_{\tau_n + \delta_n} - X^{\mu^n}_{\tau_n} | \geq \epsilon) \underset{n \rightarrow + \infty}{\longrightarrow} 0.$$ For $A>0$ which will be chosen latter, we set $$ T_A^n := \inf \left\{ t \leq T, \, |X^{\mu^n}_t| \geq A\right\}.$$ Markov's inequality yields \begin{align}\label{thmexistencemckvbisproof4}
			\P (|X^{\mu^n}_{\tau_n + \delta_n} - X^{\mu^n}_{\tau_n} | \geq \epsilon) & \leq \P ( T^n_A \leq T) + \frac{1}{\epsilon^{\beta}} \E ( |X^{\mu^n}_{\tau_n + \delta_n} - X^{\mu^n}_{\tau_n}|^\beta \1_{T_A^n \geq T}).
		\end{align}
		Notice first that \begin{align}\label{thmexistencemckvbisproof5}
			\notag\P ( T^n_A \leq T) & \leq \P ( \sup_{t \leq T} |X^{\mu^n}_t|^{\beta} \geq A^\beta) \\ &\leq \frac{R}{A^\beta}.
		\end{align}
		Then, by the triangle inequality, we have \begin{align*}
			\E ( |X^{\mu^n}_{\tau_n + \delta_n} - X^{\mu^n}_{\tau_n}|^\beta \1_{T_A^n \geq T}) & \leq \E\left( \left\vert\int_{\tau_n}^{\tau_n+ \delta_n} b_s(X^{\mu^n}_s,\mu^n_s) \, ds\right\vert^{\beta} \1_{T_A^n \geq T}\right) \\ &\quad + \E\left( \left\vert\int_{\tau_n}^{\tau_n+ \delta_n}\int_{B_1} \sigma_s(X^{\mu^n}_{s^-},\mu^n_s)z\, \NNN(ds,dz)\right\vert^{\beta} \1_{T_A^n \geq T}\right) \\&\quad + \E\left( \left\vert\int_{\tau_n}^{\tau_n+ \delta_n}\int_{B_1^c} \sigma_s(X^{\mu^n}_{s^-},\mu^n_s)z\, \NN(ds,dz)\right\vert^{\beta} \1_{T_A^n \geq T}\right) \\ &=: I_1 + I_2 + I_3
		\end{align*}
		We now estimate $I_1,$ 	$I_2$ and $I_3.$ Using the linear growth assumption on $b,$ we have for a constant $C>0$ independent of $n$ \begin{align*}
			I_1 &\leq \E\left\vert\int_{\tau_n\wedge T_A^n}^{(\tau_n+ \delta_n) \wedge T^n_A} C(1 +|X^{\mu^n}_s| + M_{\beta}(\mu^n_s)) \, ds\right\vert^\beta\\ &\leq C|\delta_n|^\beta (1 + A^\beta + R^\beta ).
		\end{align*}
		Thanks to BDG's and Jensen's inequalities, we obtain that \begin{align*}
			I_2 &\leq \E\left( \left\vert\int_{\tau_n\wedge T^n_A}^{(\tau_n+ \delta_n)\wedge T^n_A}\int_{B_1} \sigma_s(X^{\mu^n}_{s^-},\mu^n_s)z\, \NNN(ds,dz)\right\vert^{\beta} \right)\\ &\leq C \E\left( \left\vert\int_{\tau_n\wedge T^n_A}^{(\tau_n+ \delta_n)\wedge T^n_A}\int_{B_1} |\sigma_s(X^{\mu^n}_{s^-},\mu^n_s)|^2|z|^2\, \NN(ds,dz)\right\vert^{\frac{\beta}{2}} \right) \\ &\leq C(1+A^2 + R^2)^{\frac{\beta}{2}} \left(\E \int_{\tau_n\wedge T^n_A}^{(\tau_n+ \delta_n)\wedge T^n_A}\int_{B_1} |z|^2\, d\nu(z) \, ds\right)^{\frac{\beta}{2}}  \\ &\leq C(1+A^{\beta} + R^{\beta}) |\delta_n|^{\frac{\beta}{2}}.
		\end{align*}
		Since $\beta <1,$ the subadditivity of the map $\vert \cdot \vert^\beta$ yields
		
		\begin{align*}
			I_3 &\leq \E\left( \left\vert\int_{\tau_n\wedge T^n_A}^{(\tau_n+ \delta_n)\wedge T^n_A}\int_{B_1^c} \sigma_s(X^{\mu^n}_{s^-},\mu^n_s)z\, \NN(ds,dz)\right\vert^{\beta} \right)\\ &\leq  \E\left( \int_{\tau_n\wedge T^n_A}^{(\tau_n+ \delta_n)\wedge T^n_A}\int_{B_1^c} |\sigma_s(X^{\mu^n}_{s^-},\mu^n_s)|^\beta|z|^\beta\, \NN(ds,dz) \right) \\ &\leq C(1+A^\beta + R^\beta)\left(\E \int_{\tau_n\wedge T^n_A}^{(\tau_n+ \delta_n)\wedge T^n_A}\int_{B_1^c} |z|^\beta\, d\nu(z) \, ds\right)  \\ &\leq C(1+A^{\beta} + R^{\beta}) |\delta_n|.
		\end{align*}
		Using \eqref{thmexistencemckvbisproof4}, \eqref{thmexistencemckvbisproof5}, and the upper-bounds obtained previously for $I_1,$ $I_2$ and $I_3,$ we deduce that $$ 	\P (|X^{\mu^n}_{\tau_n + \delta_n} - X^{\mu^n}_{\tau_n} | \geq \epsilon) \leq \frac{R}{A^{\beta}} + \frac{C}{\epsilon^\beta}(1 + A^\beta + R^\beta) ( |\delta_n| +|\delta_n|^\beta+ |\delta_n|^{\frac{\beta}{2}}).$$ Since $R$ is fixed, we can choose $A$ large enough and then let $n$  tend to $+\infty$ to obtain that $(X^{\mu^n}_{\tau_n + \delta_n} - X^{\mu^n}_{\tau_n})_n$ converges in probability to $0.$ Thus, $([(X^{\mu^n}_t)_{t\in[0,T]}])_n$ is relatively compact in $\Pc(D_T).$ The relative compactness in $\Pc_{\beta}(D_T)$ follows from the fact that $$ \sup_{\mu \in \overline{\mathcal{B}}_R} \E \sup_{t\leq T} |X^\mu_t|^{\beta + \delta} < + \infty.$$ Indeed, this is a consequence of \cite[Proposition $2$]{fournier_moment}, since for all $t \in [0,T],$ $\mu \in \overline{\mathcal{B}}_R$ and $x\in\R^d$ $$ |b_t(x,\mu_t)| + |\sigma_t(x,\mu_t)| \leq C(1 + |x| + R),$$ and $$ \int_{B_1^c} |z|^{\beta + \delta} \, d\nu(z) <+ \infty.$$ We have proved that $\phi(\overline{\mathcal{B}}_R)$ is relatively compact in $\Pc_{\beta}(D_T).$ Thus, Schauder's fixed point theorem yields the existence of a solution to \eqref{Aedsmckv}. The moment estimate \eqref{thmexistencemckvbismoment} directly follows from \cite[Proposition $2$]{fournier_moment}.
	\end{proof}

\end{document}